\documentclass[12pt,a4paper,twoside]{article}

\newcommand{\pageformat}[6]{\setlength{\hoffset}{-1in}
                  \setlength{\voffset}{-1in}
                  \addtolength{\hoffset}{#5}
                            \addtolength{\voffset}{#6}
                            \setlength{\oddsidemargin}{#1}
                            \setlength{\evensidemargin}{#2}
                            \setlength{\textwidth}{\paperwidth}
                  \addtolength{\textwidth}{-\oddsidemargin}
                  \addtolength{\textwidth}{-\evensidemargin}
                  \addtolength{\textwidth}{-\marginparsep}
                  \addtolength{\textwidth}{-\marginparwidth}
                            \setlength{\topmargin}{#3}
                            \setlength{\textheight}{\paperheight}
                  \addtolength{\textheight}{-\topmargin}
                  \addtolength{\textheight}{-\headheight}
                  \addtolength{\textheight}{-\headsep}
                  \addtolength{\textheight}{-\footskip}
                  \addtolength{\textheight}{-#4}}
\pageformat{2cm}{3cm}{25mm}{25mm}{1pt}{0pt}

\usepackage{ifthen}
\newboolean{article}
    \setboolean{article}{true}
\newboolean{report}
\newboolean{book}
\newboolean{letter}
\newboolean{german}
\newboolean{italian}
\newboolean{nobaselinestretch}
\newboolean{nosectionappendix}
\newboolean{oldtoc}
\newboolean{nosectionequn}
\newboolean{notheorem}

\ifthenelse{\boolean{german}}{
    \usepackage{german}}{}

\usepackage[latin1]{inputenc}

\usepackage{amsmath}
\usepackage{amssymb}
\usepackage[mathscr]{eucal}

\ifthenelse{\boolean{notheorem}}{}{
    \usepackage{theorem}}

\ifthenelse{\boolean{nobaselinestretch}}{}{
    \renewcommand{\baselinestretch}{1.25}}

\newenvironment{env}[2]{\begin{#1}#2\end{#1}}{}
    \newcommand{\beq}[1]{\begin{env}{equation}{#1}}
    \newcommand{\beqn}[1]{\begin{env}{equation*}{#1}}
    \newcommand{\bal}[1]{\begin{env}{align}{#1}}
    \newcommand{\baln}[1]{\begin{env}{align*}{#1}}
    \newcommand{\bga}[1]{\begin{env}{gather}{#1}}
    \newcommand{\bgan}[1]{\begin{env}{gather*}{#1}}
    \newcommand{\bflal}[1]{\begin{env}{flalign}{#1}}
    \newcommand{\bflaln}[1]{\begin{env}{flalign*}{#1}}
    \newcommand{\bmu}[1]{\begin{env}{multline}{#1}}
    \newcommand{\bmun}[1]{\begin{env}{multline*}{#1}}
    \newcommand{\bsp}[1]{\begin{env}{split}{#1}}

    \newcommand{\eeq}{\end{env}}
    \newcommand{\eeqn}{\end{env}}
    \newcommand{\eal}{\end{env}}
    \newcommand{\ealn}{\end{env}}
    \newcommand{\ega}{\end{env}}
    \newcommand{\egan}{\end{env}}
    \newcommand{\eflal}{\end{env}}
    \newcommand{\eflaln}{\end{env}}
    \newcommand{\emu}{\end{env}}
    \newcommand{\emun}{\end{env}}
    \newcommand{\esp}{\end{env}}

\newcommand{\lf}{\vspace{2ex}}

\renewcommand{\bf}[1]{\textbf{#1}}
\renewcommand{\it}[1]{\textit{#1}}

\renewcommand{\sf}[1]{\textsf{#1}}

\newcommand{\hl}[1]{\bf{\it{#1}}}

\newcommand{\msf}[1]{\text{\small$\sf{#1}$}}

\newcommand{\eus}[1]{\mathscr{#1}}

\newcommand{\bb}[1]{\mathbb{#1}}

\newcommand{\nbd}[1]{$#1$\nobreakdash--}

\newcommand{\wt}[1]{\widetilde{#1}}

\newcommand{\vt}{\vartheta}

\newcommand{\family}[1]{\left(#1\right)}

\newcommand{\bfam}[1]{\bigl(#1\bigr)}

\newcommand{\AB}[1]{\langle#1\rangle}

\newcommand{\CB}[1]{\{#1\}}
\newcommand{\bCB}[1]{\bigl\{#1\bigr\}}
\newcommand{\BCB}[1]{\Bigl\{#1\Bigr\}}

\newcommand{\LO}[1]{(#1]}

\newcommand{\set}[2][]{
    \ifthenelse{\equal{#1}{}}{
        \CB{#2}}{
        \CB{#1~|~#2}}}
\newcommand{\bset}[2][]{
    \ifthenelse{\equal{#1}{}}{
        \bCB{#2}}{
        \bCB{#1~|~#2}}}
\newcommand{\Bset}[2][]{
    \ifthenelse{\equal{#1}{}}{
        \BCB{#2}}{
        \BCB{#1~\big|~#2}}}
\newcommand{\zero}{\CB{0}}

\DeclareMathOperator*{\limind}{lim\,ind}

\DeclareMathOperator{\id}{\normalfont\msf{id}}

\newcommand{\N}{\bb{N}}

\newcommand{\R}{\bb{R}}

\newcommand{\sB}{\eus{B}}

\newcommand{\sN}{\eus{N}}

\newcommand{\sS}{\eus{S}}

\ifthenelse{\boolean{nosectionequn}}{}{
    \numberwithin{equation}{section}
    }

\ifthenelse{\boolean{article}\or\boolean{letter}\or\boolean{nosectionequn}}{
    \setboolean{nosectionappendix}{true}}{}
\ifthenelse{\boolean{nosectionappendix}}{}{
    \renewcommand{\appendix}{
        \chapter*{\appendixname}
        \addcontentsline{toc}{chapter}{\appendixname}
        \renewcommand{\thesection}{\Alph{section}}
        \setcounter{section}{0}}}
   
\ifthenelse{\boolean{report}\or\boolean{book}}{
    }{}

\ifthenelse{\boolean{notheorem}}{}{
        \newcommand{\mnname}{Mathematical note.}
        \newcommand{\enname}{End of the note.}
        \newcommand{\definame}{Definition.}
        \newcommand{\propname}{Proposition.}
        \newcommand{\lemname}{Lemma.}
        \newcommand{\exname}{Example.}
        \newcommand{\exername}{Exercise.}
        \newcommand{\remname}{Remark.}
        \newcommand{\obname}{Observation.}
        \newcommand{\thmname}{Theorem.}
        \newcommand{\corname}{Corollary.}
        \newcommand{\proofname}{Proof.}
    \ifthenelse{\boolean{german}}{
        \renewcommand{\mnname}{Mathematische Notiz.}
        \renewcommand{\enname}{Ende der Notiz.}
        \renewcommand{\exname}{Beispiel.}
        \renewcommand{\exername}{Übung.}
        \renewcommand{\remname}{Bemerkung.}
        \renewcommand{\obname}{Beobachtung.}
        \renewcommand{\thmname}{Satz.}
        \renewcommand{\corname}{Korollar.}
        \renewcommand{\proofname}{Beweis.}}{}
    \ifthenelse{\boolean{italian}}{
        \renewcommand{\mnname}{Nota matematica.}
        \renewcommand{\enname}{Fina della nota.}
        \renewcommand{\definame}{Definizione.}
        \renewcommand{\propname}{Proposizione.}
        \renewcommand{\exname}{Esempio.}
        \renewcommand{\exername}{Esercizio.}
        \renewcommand{\remname}{Nota.}
        \renewcommand{\obname}{Osservazione.}
        \renewcommand{\thmname}{Teorema.}
        \renewcommand{\corname}{Corollario.}
        \renewcommand{\proofname}{Dimostrazione.}

       \renewcommand{\appendixname}{Appendice}

       }{}
    \theoremheaderfont{\normalfont\bfseries}
    \theoremstyle{change}
        \theorembodyfont{\rmfamily}
            \newtheorem{emp}{}[section]
                \newcommand{\bemp}[1][]{
                    \begin{emp}\hskip-\labelsep\bf{#1}\hskip\labelsep}
                \newcommand{\eemp}{\end{emp}}
\newtheorem{itemp}[emp]{}
                \newcommand{\bitemp}[1][]{
                    \begin{itemp}\hskip-\labelsep\bf{#1}\hskip\labelsep\normalfont\itshape}
                \newcommand{\eitemp}{\end{itemp}}
            \newtheorem{mn}[emp]{\mnname}
                \newcommand{\bnm}{\begin{mn}~\begin{quotation}\renewcommand{\baselinestretch}{1}\small\noindent\ignorespaces}
                \newcommand{\enm}{\end{quotation}\hfill\bf{\enname}\end{mn}}
            \newtheorem{ex}[emp]{\exname}
                \newcommand{\bex}{\begin{ex}}
                \newcommand{\eex}{\end{ex}}
            \newtheorem{exer}[emp]{\exername}
                \newcommand{\bexer}{\begin{exer}}
                \newcommand{\eexer}{\end{exer}}
            \newtheorem{defi}[emp]{\definame}
                \newcommand{\bdefi}{\begin{defi}}
                \newcommand{\edefi}{\end{defi}}
            \newtheorem{rem}[emp]{\remname}
                \newcommand{\brem}{\begin{rem}}
                \newcommand{\erem}{\end{rem}}
            \newtheorem{ob}[emp]{\obname}
                \newcommand{\bob}{\begin{ob}}
                \newcommand{\eob}{\end{ob}}
        \theorembodyfont{\normalfont\itshape}
            \newtheorem{thm}[emp]{\thmname}
                \newcommand{\bthm}{\begin{thm}}
                \newcommand{\ethm}{\end{thm}}
            \newtheorem{prop}[emp]{\propname}
                \newcommand{\bprop}{\begin{prop}}
                \newcommand{\eprop}{\end{prop}}
            \newtheorem{cor}[emp]{\corname}
                \newcommand{\bcor}{\begin{cor}}
                \newcommand{\ecor}{\end{cor}}
            \newtheorem{lem}[emp]{\lemname}
                \newcommand{\blem}{\begin{lem}}
                \newcommand{\elem}{\end{lem}}
\newenvironment{empn}[1]{\lf\noindent\bf{#1}\ignorespaces\hskip\labelsep}{\lf}
		\newcommand{\bempn}[1]{\begin{empn}{#1}}
		\newcommand{\eempn}{\end{empn}}
		\newcommand{\bitempn}[1]{\begin{empn}{#1}\normalfont\itshape}
		\newcommand{\eitempn}{\end{empn}}
                \newcommand{\bnmn}{\begin{empn}{\mnname}~\begin{quotation}\renewcommand{\baselinestretch}{1}\small\noindent\ignorespaces}
                \newcommand{\enmn}{\end{quotation}\hfill\bf{\enname}\end{empn}}
		\newcommand{\bexn}{\begin{empn}{\exname}}
		\newcommand{\eexn}{\end{empn}}
		\newcommand{\bexern}{\begin{empn}{\exername}}
		\newcommand{\eexern}{\end{empn}}
		\newcommand{\bdefin}{\begin{empn}{\definame}}
		\newcommand{\edefin}{\end{empn}}
		\newcommand{\bremn}{\begin{empn}{\remname}}
		\newcommand{\eremn}{\end{empn}}
		\newcommand{\bobn}{\begin{empn}{\obname}}
		\newcommand{\eobn}{\end{empn}}

\newcommand{\qedsymbol}{~\rule[-0.35mm]{2mm}{2mm}}
    \newcounter{proof}[emp]
    \newenvironment{Proof}[1]{
        \vspace{1ex}
        \renewcommand{\item}[1][\stepcounter{proof}(\roman{proof})]%
            {##1\hskip\labelsep}
        \noindent\textsc{#1\hskip\labelsep}}{
        \nolinebreak\qedsymbol}
    \newcommand{\proof}[1][\proofname]{
        \begin{Proof}{#1}\ignorespaces}
    \newcommand{\qed}{\end{Proof}}
    \newcommand{\noqed}{
        \renewcommand{\qedsymbol}{}
        \end{Proof}}}
    \ifthenelse{\boolean{italian}}{
        \renewcommand{\proofname}{Dimostrazione.}}{}

\usepackage[varg]{txfonts}

\begin{document}

\title{Existence of $E_0$--Semigroups for Arveson Systems: Making Two Proofs into One}
\author{}
\author{
~\\
Michael Skeide
\\
}
\date{May 2006}

{
\renewcommand{\baselinestretch}{1}
\maketitle



\begin{abstract}
\noindent
Since quite a time there were available only two rather difficult and involved proofs, the original one by Arveson and a more recent one by Liebscher, of the fact that for every Arveson system there exists an \nbd{E_0}semigroup. We put together two recent short proofs, one by Skeide and one by Arveson, to obtain a still simpler one, which unfies the advantages of each proof and discards with their disadvantages.
\end{abstract}

}


\section{Introduction}

In \cite{Arv89} Arveson associates with every \nbd{E_0}semigroup $\vt$ (that is, a strongly continuous semigroup of unital normal endomorphisms of the algebra $\sB(H)$ of all adjointable operators on a separable infinite-dimensional Hilbert space $H$) a product system $E^\otimes=\bfam{E_t}_{t\in(0,\infty)}$ of Hilbert spaces $E_t$, a so-called \it{Arveson system}. As Arveson systems classify \nbd{E_0}semigroups up to cocycle conjugacy it is natural to ask whether every Arveson system may be obtained as the Arveson system of an \nbd{E_0}semigroup. Arveson answered this question in the affirmative sense in \cite{Arv90}. However, the proof is long and involves deep analytic techniques, some of which first had to be developed in \cite{Arv90a,Arv89a}. Also a second proof due to Liebscher \cite{Lie00p1} appears to be involved.

Recently, Skeide \cite{Ske06p2} and, shortly after, Arveson \cite{Arv06p} have given short proofs of this result. The idea of the proof in \cite{Ske06p2} is plain and unitality of the constructed endomorphisms is obvious, while the verification of the semigroup property is is rather tedious. The proof in \cite{Arv06p} has no problems with the semigroup property, while the verification of unitality requires a computation. In these notes we show that the constructions from \cite{Ske06p2} and from \cite{Arv06p} actually are unitarily equivalent. In this way, we can avoid in each proof that part which is less obvious from its construction.

The main accent in this short note is on establishing unitary equivalence of the two constructions. For this to us it appears more convenient to discuss first \cite{Ske06p2} and then switch to \cite{Arv06p}. Of course, the whole thing could be prepared also in the opposite direction. We note also that here we describe the construction from \cite{Ske06p2} with all orders in tensor products reversed (construction of a right dilation instead of a left dilation). This operation does not cause any complication but facilitates then comparison with the construction in \cite{Arv06p}. We should also say that we leave all details about verification of measurabilities that go beyond the necessary conditions on square integrability (that is in particular mesaurability) of certain sections to either of the articles \cite{Ske06p2,Arv06p}; see also Remark \ref{mrem}. Our emphasis is on algebraic problems like associativity and unitality.

\section{Preliminaries}

Throughout these notes we assume that $E^\otimes=\bfam{E_t}_{t\in(0,\infty)}$ is a fixed \hl{Arveson system} of Hilbert spaces $E_t$. Algebraically, this means that we have associative identifications $u_{t,s}\colon E_t\otimes E_s\rightarrow E_{t+s}$. Technically, we assume that the bundle $E^\otimes$ has a Borel structure isomorphic to the trivial bundle $(0,\infty)\times H_0$ for an infinite-dimensional separable Hilbert space $H_0$, and we assume that for every pair of measurable functions $x\colon t\mapsto x(t)\in E_t\subset E^\otimes$ and $y\colon t\mapsto y(t)\in E_t\subset E^\otimes$ also the function $(t,s)\mapsto u_{t,s}(x(t)\otimes y(s))\in E_{t+s}\subset E^\otimes$ is measurable. Every Arveson system in the sense of \cite[Definition 1.4]{Arv89} fulfills this condition. We will follow the conventions in \cite{Arv89} and write $x(t)y(s)$ for $u_{t,s}(x(t)\otimes y(s))$. Associativity simply means then that this \it{product} on $E^\otimes$ is associative.

By \cite{Ske06p2} a \hl{right dilation} of an Arveson system $E^\otimes$ is a nontrivial (and, therefore, infinite-dimensional) separable Hilbert space $R$ with unitary identifications $w_t\colon E_t\otimes R\rightarrow R$ which behave associatively with respect to the product system structure. Writing $xh$ for $w_t(x\otimes h)$, this means just that $(xy)h=x(yh)$ for all $x\in E_t,y\in E_s,h\in R$. In other words, $R$ is a left module over the ring generated by the semigroup $E^\otimes$.

\brem\label{E_0}
A right dilation of $E^\otimes$ induces a semigroup $\vt=\bfam{\vt_t}_{t\in\R_+}$ of normal unital endomorphism $\vt_t(a):=w_t(\id_t\otimes a)w_t^*$ (and $\vt_0=\id_{\sB(R)}$) of $\sB(R)$. It is easy to see that the Arveson system associated to $\vt$ as in \cite{Arv89} is $E^\otimes$.
\erem

\brem
A right dilation induces an essential (that is, nondegenerate) representation $\Phi=\bfam{\Phi_t}_{t\in(0,1)}$ of $E^\otimes$ on $R$, namely, $\Phi_t(x)h:=xh$. Conversely, if $\Phi$ is an essential representation of $E^\otimes$ on $R\ne\zero$ (separable), then $w_t\colon x\otimes h\mapsto\Phi_t(x)h$ defines a right dilation. By Remark \ref{E_0} an essential representation of an Arveson system induces, therefore, also a unital endomorphism semigroup having this Arveson system.
\erem

\brem\label{mrem}
So far we spoke about right dilations in algebraic terms. For that the endomorphism semigroup induced by a right dilation be strongly continuous, that is, be an \hl{\nbd{E_0}semi\-group}, it is sufficient that the right dilation be \hl{measurable} in the sense that for every pair of measurable functions $x\colon t\mapsto x(t)\in E_t$ and $h\colon t\mapsto h(t)\in R$ the function $t\mapsto x(t)h(t)$ be measurable. See \cite{Ske06p2,Arv06p} for possibilities for how to conclude from this to continuity of the semigroup. In \cite[Proposition 4.1]{Ske06p2} there is a self-contained proof (using only the fact that unitary groups on a separable Hilbert space are strongly continuous if they are weakly measurable) of that a measurable (left or right) dilation gives rise to an \nbd{E_0}semigroup. A similar result is \cite[Proposition 2.7]{Arv89}.
\erem

\section{The first construction}\label{1st}

The construction in \cite{Ske06p2} follows two steps (here rephrased suitably in terms of right dilations). In the first step, one constructs a right dilation of the discrete subsystem $\bfam{E_n}_{n\in\N}$ of $E^\otimes$. In other words, one constructs a separable Hilbert space $\breve{R}\ne\zero$ and unitary identifications $\breve{w}_n\colon E_n\otimes\breve{R}\rightarrow\breve{R}$ that compose associatively with the product system structure. The construction of such dilations for discrete product systems is easy and we come back to it in a minute, because for comparison with \cite{Arv06p} we need to make a concrete choice. Having the right dilation of the discrete subsystem, the idea of \cite{Ske06p2} is to put $R:=\bfam{\int_0^1E_\alpha\,d\alpha}\otimes\breve{R}$ and to write down for every $t\in(0,\infty),n:=\CB{t}$ (the unique integer such that $t-n\in\LO{0,1}$) specific versions of the following isomorphisms
\bal{\notag
E_t\otimes R
~=~
E_t\otimes\family{\int_0^1E_\alpha\,d\alpha}\otimes\breve{R}
&
~\cong~
\family{\int_t^{1+t}E_\alpha\,d\alpha}\otimes\breve{R}
\\[2ex]\notag
&
~\cong~
\family{\int_{t-n}^1E_\alpha\,d\alpha}\otimes E_n\otimes\breve{R}
~~\oplus~~
\family{\int_0^{t-n}E_\alpha\,d\alpha}\otimes E_{n+1}\otimes\breve{R}
\\[2ex]\label{idea}
&
~\cong~
\family{\int_{t-n}^1E_\alpha\,d\alpha}\otimes\breve{R}
~~\oplus~~
\family{\int_0^{t-n}E_\alpha\,d\alpha}\otimes\breve{R}
~~\cong~~
R
}\eal
and show that they iterate associatively. This program works for an arbitrary right dilation of the discrete subsystem, not only for the one we shall consider in the sequel, which leads to something unitarily equivalent to \cite{Arv06p}. But the verification of associativity is tedious; see \cite{Ske06p2}.

The concrete right dilation of the discrete subsystem suggested in \cite{Ske06p2} is obtained in the following way. Choose a unit vector $e$ in $E_1$. For $n,m\in\N$ define isometric embeddings $E_m\rightarrow E_{m+n}$ by $x\mapsto xe^n$. (That is, we identify $E_m$ as the subspace $E_me^n=u_{m,n}(E_m\otimes e^n)$ of $E_{m+n}$.) These embeddings form an inductive system. Let us denote by $\breve{R}$ the (completed) inductive limit. It is easy to check that the family $\bfam{u_{m,n}}_{m\in\N}$ is compatible with the inductive limit over $m\in\N$, that it defines a unitary $\breve{w}_n\colon E_n\otimes\breve{R}\rightarrow\breve{R}$ and that the family of all these $\breve{w}_n$ is a right dilation.

So far the construction from \cite{Ske06p2}. For the following sections it is important to observe that in the construction of $R$ nobody prevents us from exchanging the order of inductive limit and direct integral. So let us define the spaces $K_m:=\int_m^{m+1}E_\alpha\,d\alpha\cong\bfam{\int_0^1E_\alpha\,d\alpha}\otimes E_m=:\breve{K}_m$. The isometric embeddings $E_m\rightarrow E_me^n\subset E_{m+n}$ give rise to isometric embeddings, first, of $\breve{K}_m\rightarrow\breve{K}_{m+n}$ and, then, of $K_m\rightarrow K_{m+n}$. Clearly, the inductive limit $\breve{K}:=\limind_m\breve{K}$ over the spaces $\breve{K}_m$ is just $R$ as we used it above. And because the canonical identifications $\bfam{f(\alpha)}_{\alpha\in\LO{0,1}}\otimes x_m$ and $\bfam{f(\alpha-m)x_m}_{\alpha\in\LO{m,m+1}}$ of $\breve{K}_m$ and $K_m$ are compatible with the inductive structures, also $R$ and the inductive limit $K:=\limind_m K_m$ over the spaces $K_m$ are canonically isomorphic. We will investigate this latter inductive limit $K$ in the Section \ref{CSEC}. In particular, we will see that it coincides with the Hilbert space constructed in \cite{Arv06p}.

We close this section by analyzing how the dilation constructed in \cite{Ske06p2} looks like in terms of these inductive limits. A typical element of $R=\bfam{\int_0^1E_\alpha\,d\alpha}\otimes\breve{R}$ is $\bfam{f(\alpha)}_{\alpha\in\LO{0,1}}\otimes y$ with $f(\alpha)\in E_\alpha$ and $y\in\breve{R}$. The unitaries $w_t$ from \cite{Ske06p2} as suggested by Equation \eqref{idea} act on $x\otimes\bfam{f(\alpha)}_{\alpha\in\LO{0,1}}\otimes y$ in $E_t\otimes R$ pointwise on sections. Putting $n=\CB{t+\alpha}$, the unitary $w_t$ sends the point $x\otimes f(\alpha)\otimes y$ $(\alpha\in\LO{0,1})$ in the section to
\beqn{
(\id_{t+\alpha-n}\otimes\breve{w}_n)(u_{t+\alpha-n,n}^*(xf(\alpha))\otimes y)
~\in~
E_{t+\alpha-n}\otimes\breve{R}.
}\eeqn
That is, $w_t$ sends $x\otimes\bfam{f(\alpha)}_{\alpha\in\LO{0,1}}\otimes y$ to a section in $R$ which at time $\beta=t+\alpha-n\in\LO{0,1}$ assumes the value $(\id_{t+\alpha-n}\otimes\breve{w}_n)(u_{t+\alpha-n,n}^*(xf(\alpha))\otimes y)$.

Now recall that the inductive limit $\breve{R}$ is generated by the increasing sequence of subspaces $E_m$. If $y\in E_m$, then $\breve{w}_n$ sends the tensor product of an element in $E_n$ with $y$ to an element in $E_{n+m}\subset\breve{R}$. So on the level of the spaces $\breve{K}_m$ the point $x\otimes f(\alpha)\otimes y$ of the section $x\otimes\bfam{f(\alpha)}_{\alpha\in\LO{0,1}}\otimes y\in E_t\otimes\breve{K}_m$ ends up in a point of a section in $\breve{K}_{n+m}$. On the level of the spaces $K_m$ the point $x\otimes f(\alpha-m)y$  $(\alpha\in\LO{m,m+1})$ of a section in $E_t\otimes K_m$ ends up on the point $xf(\alpha-m)y$ at $\beta=t+\alpha\in\LO{n+m,n+m+1}$ of a section in $K_{n+m}$ where $n=\CB{t+\alpha-m}$. In other words, if now $\bfam{f(\alpha)}_{\alpha\in\LO{m,m+1}}$ is an arbitrary section in $K_m$, then the point $x\otimes f(\alpha)$ ends up in the point $xf(\alpha)$ of a section in $K_{n+m}$ at $\beta$ with $\beta$ and $n$ as before.

\section{The second construction}\label{2nd}

In \cite{Arv06p} Arveson constructs a Hilbert space $H$ as follows. Let $\sS$ denote the space of all locally square integrable sections $f=\bfam{f(t)}_{t\in(0,\infty)}\subset E^\otimes$ which are \hl{stable} with respect to the unit vector $e\in E_1$, that is, for which there exists an $\alpha_0>0$ such that
\beqn{
f(\alpha+1)
~=~
f(\alpha)e
}\eeqn
for all $\alpha\ge\alpha_0$. By $\sN$ we denote the subspace of all sections in $\sS$ which are eventually $0$, that is, of all sections $f\in\sS$ for which there exists an $\alpha_0>0$ such that $f(\alpha)=0$ for almost all $\alpha\ge\alpha_0$. A straightforward verification shows that
\beqn{
\AB{f,g}
~:=~
\lim_{m\to\infty}\int_m^{m+1}\AB{f(\alpha),g(\alpha)}\,d\alpha
}\eeqn
defines a semiinner product on $\sS$ and that $\AB{f,f}=0$ if and only if $f\in\sN$. (Actually, we have
\beqn{
\AB{f,g}
~=~
\int_T^{T+1}\AB{f(\alpha),g(\alpha)}\,d\alpha
}\eeqn
for all sufficiently large $T>0$; see \cite[Lemma 2.1]{Arv06p}.) So, $\sS/\sN$ becomes a pre-Hilbert space with inner product $\AB{f+\sN,g+\sN}:=\AB{f,g}$. By $H$ we denote its completion.

After these preparations it is completely plain to see that for every $t>0$ the map $x\otimes f\mapsto xf$, where
\beqn{
(xf)(\alpha)
~=~
\begin{cases}
xf(\alpha-t)&\alpha>t,
\\
0&\text{else},
\end{cases}
}\eeqn
defines an isometry $E_t\otimes H\rightarrow H$, and that these isometries iterate associatively. Surjectivity of these isometries is slightly less obvious.

\section{Comparison and integration of the two approaches}\label{CSEC}

We claim that $R$ from Section \ref{1st} and $H$ from Section \ref{2nd} are canonically isomorphic in a way such that the mappings $E_t\otimes R\rightarrow R$ and $E_t\otimes H\rightarrow H$ become unitarily equivalent. This shows immediately that the former iterate associatively (because the latter do) and that the latter are unitaries (because the former are). In this way, we remove from each construction the most tedious verifications.

We have already established in the end of Section \ref{1st} that the space $R$ can be viewed as inductive limit $\breve{K}$ over the spaces $\breve{K}_m$, that this inductive limit is canonically isomorphic to the inductive limit $K$ over the spaces $K_m=\int_m^{m+1}E_\alpha\,d\alpha$ and we have established how the $w_t$ act on $E_t\otimes R$ when restricted to $E_t\otimes K_m$.

Le $f=\bfam{f(\alpha)}_{\alpha\in\LO{m,m+1}}$ be in $K_m$ and define the section $\wt{f}\in\sS$ by setting
\beqn{
\wt{f}(\alpha)
~=~
\begin{cases}
0&\alpha\le m,
\\
f(\alpha-n)e^n&n\in\N_0,m+n<\alpha\le m+n+1.
\end{cases}
}\eeqn
Then $f\mapsto\wt{f}+\sN$ defines an isometry $K_m\rightarrow H$. Morevoer, recalling that the inductive structure of the family of spaces $K_m$ is given by embeddings that embed the section $f\in K_m$ as the section $\bfam{f(\alpha-n)e^n}_{\alpha\in\LO{m+n,m+n+1}}$ into $K_{n+m}$, we easily check that $K_m\rightarrow H$ and $K_m\rightarrow K_{m+n}\rightarrow H$ coincide. In other words, we have an isometry from $R\cong K$ into $H$. Moreover, if $\wt{f}$ is in $\sS$, then there is an $\alpha_0>0$ so that $\wt{f}(\alpha+1)=\wt{f}(\alpha)e$ for all $\alpha\ge\alpha_0$. In other words, if we choose an integer $m\ge\alpha_0$, then up to an element in $\sN$ the section $\wt{f}$ is the image of the section $f\in K_m$ defined by setting $f(\alpha)=\wt{f}(\alpha)$ $(\alpha\in\LO{m,m+1})$. In other words, by the imbedding of $R$ into $H$ we obtain a total subset of $H$. Therefore, actually we have defined a unitary $R\rightarrow H$. (Note that, actually, we have identified $\bigcup_{m\in\N}K_m$ with $\sS/\sN$.)

Now recall that $w_t$ sends $x\otimes f$ with a section $f$ in $\breve{K}_m$ or, what is the same up to canonical isomorphism, in $K_m$ to a section that, for some $n\in\N_0$, lies partly in $K_{n+m}$ partly in $K_{n+1+m}$. In fact, we may split $f=\bfam{f(\alpha)}_{\alpha\in\LO{m,m+1}}$ into two parts $\bfam{f(\alpha)}_{\alpha\in\LO{m,m+t-n}}+\bfam{f(\alpha)}_{\alpha\in\LO{m+t-n,m+1}}$ with $n=\CB{t}$ so that the first part ends up in $K_{n+m}$ while the second part ends up in $K_{n+1+m}$. (Note that, actually, we have $K_{n+m}=K_{n+m}e\subset K_{n+m+1}$ and the fact that $w_t$ is isometric shows even that the two parts of $x\otimes f$ end up in orthogonal parts of $K_{n+1+m}$. But this is not the point.) We simplify life by noting that it is sufficient to consider only those sections for which one of the parts is zero. Let us denote the result in $K_{n+m}$ or in $K_{n+1+m}$ by $xf$.

It is now completely plain to verify that the sections $x\wt{f}$ and $\wt{xf}$ in $\sS$ coincide eventually and, therefore, $x\tilde{f}+\sN$ and $\wt{xf}+\sN$ coincide in $H$.

\newcommand{\Swap}[2]{#2#1}\newcommand{\Sort}[1]{}
\providecommand{\bysame}{\leavevmode\hbox to3em{\hrulefill}\thinspace}
\providecommand{\MR}{\relax\ifhmode\unskip\space\fi MR }
\providecommand{\MRhref}[2]{%
  \href{http://www.ams.org/mathscinet-getitem?mr=#1}{#2}
}
\providecommand{\href}[2]{#2}

\setlength{\baselineskip}{2.5ex}



\end{document}